\documentclass{amsart}

\usepackage{amssymb,color}
\usepackage{amsfonts}
\usepackage{amsmath}
\usepackage{euscript}
\usepackage{enumerate}
\usepackage{pdfsync}
\newtheorem*{Theorem1'}{Theorem 1'}

\theoremstyle{definition}

\theoremstyle{remark}

\begin{document}

\title[On the normalizer of an iterated wreath product]{On the normalizer of an iterated wreath product}

\author{Fernando Szechtman}
\address{Department of Mathematics and Statistics, University of Regina, Canada}
\email{fernando.szechtman@gmail.com}
\thanks{The author was partially supported by NSERC grant RGPIN-2020-04062}




\begin{abstract} Given a group $G$ and $n\geq 0$, let $W(G,n)$ be the associated iterated wreath product -- unrestricted
when $G$ is infinite -- viewed
as a permutation group on $G^n$. We prove that the normalizer of $W(G,n)$ in the symmetric group $S(G^n)$ is equal to 
$M_n\ltimes W(G,n)$, 
where $M_n$ is isomorphic to~$\mathrm{Aut}(G)^n$. The action of $\mathrm{Aut}(G)^n$ on $W(G,n)$ is recursively described.
\end{abstract}

\maketitle

\section{Introduction}

Given a group $G$ and $n\geq 0$, we let $W(G,n)$ stand for the associated iterated standard wreath product, taken in the
unrestricted sense when $G$ is infinite, where
$G$ is viewed as a permutation group on itself via the regular representation and, correspondingly, $W(G,n)$ is viewed 
as a permutation group on $G^n$.

If  $p$ is a prime number, $W(C_p,n)$ is a Sylow $p$-subgroup of the symmetric group $S_{p^n}$ \cite{K}, and the normalizer
of $W(C_p,n)$ in $S_{p^n}$ is equal to $M_n\ltimes W(C_p,n)$, where $M_n$ is isomorphic to $C_{p-1}^n$ \cite[Appendix 13]{BJ}. As it turns out, $C_p$ is but the 
simplest group for which the latter holds.

\medskip

\noindent{\bf Theorem.} The normalizer
of $W(G,n)$ in the symmetric group $S(G^n)$ is equal to $M_n\ltimes W(G,n)$, where $M_n$ is isomorphic to $\mathrm{Aut}(G)^n$.

The action of $\mathrm{Aut}(G)^n$ on $W(G,n)$ is recursively described in Section 3.

\medskip

The problem of computing normalizers in the general context of finite permutation groups is addressed in \cite{CCR} from
an algorithmic standpoint (see also references therein). The normalizers of certain finite subgroups of the so-called homogeneous symmetric groups are determined
in \cite{GK}. The automorphism group of $W(C_p,n)$ for $p$ odd is studied in \cite{B} and \cite{L}, and the 
the automorphism group of the normalizer $W(C_p,n)$ in $S_{p^n}$ is described in \cite{BL}. A variation of $W(C_p,n)$ with first
term $C_{p^e}$, as well as the corresponding automorphism group, is investigated in \cite{R}. The irreducible representations
of the iterated wreath product of finite cyclic groups are elucidated in \cite{IW}.

\section{Unrestricted wreath product of permutation groups}

Given permutation groups $H$ and $K$ acting on sets $X$ and $Y$ respectively, the unrestricted wreath product $H\wr K$ of $H$ and $K$
is defined as a permutation group on $Z=X\times Y$ as follows. 

For $y\in Y$ we set $x_y=(x,y)$, where $x\in X$, and let
$X_y=\{x_y\,|\, x\in X\}$, so that ${\mathcal P}=\{X_y\,|\, y\in Y\}$ is partition of $Z$. Given $h\in H$ and $y\in Y$ we define 
$h_y\in S(X_y)$ by $x_y h_y=(xh)_y$ for all $x\in X$ and let $H_y=\{h_y\,|\, h\in H\}$, whence
the maps $h\mapsto h_y$ and $x\mapsto x_y$ yield a similarity from $H$ to $H_y$.

We naturally view each $S(X_y)$ as a subgroup of $S(Z)$.
Let $S_Y(Z)$ be the subgroup of $S(Z)$ of all permutations $f$ of $Z$ 
such that $X_y f=X_y$ for all $y\in Y$,
and let $B$ be the subgroup of all $f\in S_Y(Z)$ such that the restriction $f_y$ of $f$ to $X_y$ is in~$H_y$ for each $y\in Y$.

We have an embedding $\Lambda:K\to S(Z)$, where $k\Lambda$ is given by $x_y\mapsto x_{yk}$ for all $x\in X$, $y\in Y$, and $k\in K$. If no risk
of confusion is possible we will identify $K$ with $K\Lambda$. Thus 
$$
S(X_y)^{k}=S(X_{yk}),\; h_y^k=h_{yk},\; H_y^k=H_{yk}, \quad h\in H, y\in Y, k\in K.
$$
It follows that $K$ normalizes $B$ with $B\cap K=1$, and we set $W=H\wr K=K\ltimes B$ and $N=N_{S(Z)} W$.
The group $B$ is known as the base group of $W$ and is isomorphic to the Cartesian product of the $H_y$, $y\in Y$.
Note that if $Y$ is finite then $B$ is the internal direct product of its subgroups $H_y$, $y\in Y$. Observe also that if
$|X|=1$ then $W$ is similar to $K$ (this case arises in the inductive proof of our main result when passing from $n=0$ to $n=1$).

\medskip

\noindent{\bf Lemma A.} If $H$ and $K$ are  transitive then so is $W$. 

\medskip

\noindent{\sc Proof.} Given $x,x'\in X$ and $y,y'$ there are $h\in H$ and $k\in K$ such that $xh=x'$ and $yk=y'$,
so $x_y h_y k=x'_{y'}$.

\medskip

\noindent{\bf Lemma B.} If $K$ is semiregular, then conjugation by any element of $N$ sends every $H_y$ back into~$B$.

\medskip

\noindent{\sc Proof.} Suppose, if possible, that there is $\alpha\in N$ such that
for some $y\in Y$  and $h\in H$  we have $h_y^\alpha\notin B$.
Then $h_y^\alpha=k f$, where $1\neq k\in K$ and $f\in B$. For $\gamma\in S(Z)$ let $M(\gamma)$ be the subset of $Z$ of points moved by $\gamma$.
Thus $M(kf)=Z$ and $M(h_y)\subseteq X_y$, so
$$
Z\alpha=Z=M(kf)=M(h_y^\alpha)=M(h_y)\alpha\subseteq  X_y\alpha,
$$
which implies $Y=\{y\}$ and a fortiori $K=1$, a contradiction.

\medskip

\noindent{\bf Lemma C.} If $K$ is semiregular and $H$ is transitive, then the partition ${\mathcal P}$ of $Z$ is $N$-stable.
Thus, there is a group homomorphism $\Omega:N\to S({\mathcal P})$, say $\alpha\mapsto\widetilde{\alpha}$, with kernel $N\cap S_Y(Z)$.
\medskip

\noindent{\sc Proof.}  Let $\alpha\in N$ and $y\in Y$. Take any $x\in X$.
Then $x_y\alpha\in X_{y'}$ for a unique $y'\in Y$. Given any other $x'\in X$, we have $x'=xh$ for some $h\in H$,
so that $x'_y=x_y h_y$. By Lemma B, $h_y^\alpha=f\in B$, so
$$
x'_y\alpha=x_y h_y \alpha=x_y \alpha f \in X_{y'}.
$$
Thus $X_y\alpha\subseteq X_{y'}$ and therefore  $X_y\subseteq X_{y'}\alpha^{-1}$. But $\alpha^{-1}\in N$
and $X_{y'}\alpha^{-1}\cap X_y\neq\emptyset$, so by above $X_{y'}\alpha^{-1}\subseteq X_{y}$. Hence $X_{y'}\alpha^{-1}=X_{y}$
and therefore $X_y\alpha=X_{y'}$.

\medskip

We suppose for the remainder of this section that $K$ is regular. 
Thus, we may assume without loss that $Y=K$ under the regular representation. 

\medskip

\noindent{\bf Lemma D.} The map $\Gamma:\mathrm{Aut}(K)\to N$ given
by
$$
x_k\gamma\Gamma=x_{k\gamma},\quad x\in X, k\in K, \gamma\in \mathrm{Aut}(K),
$$
is a group monomorphism satisfying
$$
h_k^{\gamma \Gamma}=h_{k\gamma},\quad h\in H, k\in K, \gamma\in \mathrm{Aut}(K),
$$
$$
k \Lambda^{\gamma\Gamma}=k\gamma \Lambda,\quad k\in K,\gamma\in \mathrm{Aut}(K).
$$

\medskip

\noindent{\sc Proof.} This is a routine calculation.

\medskip

We suppose for the remainder of this section that $H$ is transitive.

\medskip

\noindent{\bf Lemma E.}  The image of $\Omega:N\to S({\mathcal P})$ is 
$\widetilde{\mathrm{Aut}(K)\Gamma}\ltimes\widetilde{\Lambda K}$.

\medskip

\noindent{\sc Proof.} As $B\subseteq\ker\Omega$, we see that $\widetilde{N}$ is contained in the normalizer
of the regular subgroup $\widetilde{\Lambda K}$. It is well-known \cite[Ch. 6]{H} that this normalizer is the semidirect
product of $\widetilde{\Lambda K}$ by the group of permutations of ${\mathcal P}$ associated to the automorphisms of $K$.
This group is precisely $\widetilde{\mathrm{Aut}(K)\Gamma}$.

\medskip

If no risk of confusion is possible, we will write $\mathrm{Aut}(K)\ltimes K$ instead of $\mathrm{Aut}(K)\Gamma\ltimes\Lambda K$.

\newpage

\noindent{\bf Lemma F.}  We have $N=(\mathrm{Aut}(K)\ltimes K)(N\cap S_K(Z))$.
\medskip 

\noindent{\sc Proof.} This is a consequence of Lemmas C and E. 

\medskip   

We assume in what follows that $N_{S(X)} H =M\ltimes H$ for some subgroup $M$ of $S(X)$. To each $m\in M$ 
we associate $m^*\in N\cap S_K(Z)$ defined by $x_k m^*=(xm)_k$ for all $x\in X$ and $k\in K$. It is clear
that the map $m\mapsto m^*$ yields a group monomorphism $M\to N\cap S_K(Z)$ whose image will be denoted by $M^*$.

\medskip

\noindent{\bf Lemma G.}  We have $N=M^*\mathrm{Aut}(K)\ltimes W$, where $\mathrm{Aut}(K)\cap M^*=1$,
$[M^*,K]=1$, and $[M^*,\mathrm{Aut}(K)]=1$.

\medskip 

\noindent{\sc Proof.} Recall that for $f\in S_K(Z)$, $h\in H$, and $k\in K$, we write $f_k$ for the restriction of $f$ to $X_k$,
and $h_k\in H_k$ for the map $x_k\mapsto (xh)_k$. Given $m\in M$ and $k\in K$, we let $m_k\in S(X_k)$ be given by $x_k m_k=(xm)_k$ for 
all $x\in X$, and set $M_k=\{m_k\,|\, m\in M\}$. Thus $N_{S(X_k)} H_k=M_k\ltimes H_k$. 

We claim that $N\cap S_K(Z)=M^* B$. Clearly $M^* B\subseteq N\cap S_K(Z)$. Let $\alpha\in N\cap S_K(Z)$.
Given any $k\in K$, we have $H_k^\alpha\subseteq W\cap S(X_k)=H_k$, whence $\alpha_k\in M_k\ltimes H_k$. Thus, there
is $b\in B$ such that $\beta=\alpha b$ satisfies $\beta_k\in M_k$ for all $k\in K$. We assert that $\beta=m^*$
for some $m\in M$, so that $\alpha=m^* b^{-1}\in M^* B$, whence $N\cap S_K(Z)\subseteq M^* B$, thereby proving the claim.
To see the assertion, let $k\in K$. As $\beta_k\in M_k$, there exists $m^k\in M$ (where the superscript indicates 
dependence on $k$) such that $\beta_k=m^k_k$, that is, $x_k \beta_k=(xm^k)_k$ for all $x\in X$. Setting $m=m^1$,
we proceed to show that $\beta=m^*$. As $k\Lambda\in W$ and $\beta\in N$, we have $(k\Lambda)^\beta\in W$, so
there exist $\ell^k\in K$ and $b^k\in B$ (where the superscripts indicate
dependence on $k$) such that $(k\Lambda_k)^\beta=(\ell^k\Lambda) b^k$. Here $\beta,k\Lambda, \ell^k\Lambda, b^k\in N$,
with $\beta,b^k\in S_K(Z)$, so making use of the homomorphism from Lemma C, we infer $\ell^k=k$. On the other hand,
as $b^k\in B$, the very definition of $B$ ensures the existence of $h^k\in H$ (where the superscript indicates 
dependence on $k$) such that $b^k_k=h^k_k$, that is, $(x b^k)_k=(x h^k)_k$ for all $x\in X$. Thus for every $x\in X$, we have
$$
\begin{aligned}
(x m^k)_k &=x_k m^k_k=x_k \beta_k=x_k \beta=x_1 (k\Lambda)\beta=x_1 \beta (k\Lambda) b^k=x_1 \beta_1 (k\Lambda) b^k\\
&=(x m)_1 (k\Lambda) b^k=(xm)_k b^k=(xm)_k b^k_k=(xm)_k h^k_k=(xm h^k)_k.
\end{aligned}
$$
Hence $m^k=m h^k$ for all $k\in K$. As $m^k, m\in M$ and $h^k\in H$, with $M\cap H=1$, we deduce $m^k=m$ for all $k\in K$.
Therefore
$$
x_k \beta=x_k\beta_k=x_k m^k_k=x_k m_k=(xm)_k=x_k m^*,\quad k\in K,
$$
which implies $\beta=m^*$, as stated. It now follows from Lemma F that $N=\mathrm{Aut}(K)K M^* B$.
It is easy to see that $[M^*,K]=1$ and $[M^*,\mathrm{Aut}(K)]=1$, whence $N=M^*\mathrm{Aut}(K)W$. Here $M^*$ normalizes $W$
by construction, and $\mathrm{Aut}(K)$ normalizes $W$ by Lemma  D. 

Suppose, finally, that $m^*\gamma=k b$ for some $m\in M$, $\gamma\in \mathrm{Aut}(K)$, $k\in K$, and $b\in B$.
Applying the homomorphism from Lemma C, we see that $k=1$ and $\gamma=1$. Then $m^*=b$, so $M\cap H=1$
forces $m=1$ and $b=1$.

\medskip   

\section{Proof of the theorem and action of $\mathrm{Aut}(G)^n$ on $W(G,n)$}

\medskip 

\noindent{\bf Lemma H.} The unrestricted iterated wreath product $W(G,n)$ is a transitive subgroup of $S(G^n)$.

\medskip 

\noindent{\sc Proof.} We argue by induction on $n$. The result is clear if $n=0$.
Suppose $W(G,n)$ is a transitive subgroup of $S(G^n)$ for some $n\geq 0$. 
Set $Z=G^{n+1}=X\times Y$, where $X=G^n$ and $Y=G$, and take $H=W(G,n)$ and $K=G$. Then 
$W(G,n+1)$ is a transitive subgroup of~$S(Z)$ by Lemma A.

\medskip

\noindent{\sc Proof of the theorem.} We argue by induction on $n$, the theorem being trivially true when $n=0$.
Suppose the result is true for some $n\geq 0$. Adopting the notation of Lemma H,
the inductive hypothesis yields $N_{S(X)} W(G,n)=M\ltimes W(G,n)$,
where $M\cong\mathrm{Aut}(G)^n$. By Lemma H, $W(G,n)$ is transitive, whence Lemma G gives
$N=M^*\mathrm{Aut}(G)\ltimes W(G,n+1)$, with $M^*\cong M$, $M^*\cap\mathrm{Aut}(G)=1$, and $[M^*,\mathrm{Aut}(G)]=1$, so
$M^*\mathrm{Aut}(G)\cong \mathrm{Aut}(G)^{n+1}$.

\medskip

Let us recursively indicate the action of $\mathrm{Aut}(G)^n$ on $W(G,n)$. If $n=0$ there is nothing to~do.
Given $n\geq 0$, we have $N_{S(X)}W(G,n)=M\ltimes W(G,n)$ in the notation of Lemma H, where the action
of $M\cong \mathrm{Aut}(G)^n$ on $W(G,n)$ is assumed to be known. According to Lemma G, we have
$N_{S(Z)} W(G,n+1)=(M^*\times \mathrm{Aut}(G))\ltimes (G\ltimes B)$ with $[M^*,G]=1$. The action
of $\mathrm{Aut}(G)$ on $G\ltimes B$ is given in Lemma D, while the action of $M^*$ on $B$ is described prior to Lemma G.

\medskip

\noindent{\bf Acknowledgment.} We thank the referee for a careful reading of the manuscript.


\end{document}